\documentstyle[12pt]{article} 

\newcommand\pone{{\bf {\rm P\/}^1}}
\newcommand\ptwo{{\bf {\rm P\/}^2}}
\newcommand\pthree{{\bf {\rm P\/}^3}}
\newcommand\pfour{{\bf {\rm P\/}^4}}

\newcommand\lm{\lambda}
\newcommand{\binom}[2] 
{ \left(\begin{array}{c}#1\\#2\end{array}\right)}

\begin{document}

\title{A smooth surface in $\pfour$ not of general type has degree at most 66}
\author{Robert Braun and Michele Cook}
\date{}
\maketitle

This is a continuation of the papers of Braun and Fl{\o}ystad [BF] and 
Cook [C]  bounding the degree of smooth surfaces not of general 
type in $\pfour$.

We will prove the following: 

\vspace{3mm}
{\bf Theorem 1.}

{\em Let $S$ be a smooth surface of degree $d$ in $\pfour$ not of general type. Then $d \leq 66$. }

\vspace{3mm}
The new idea in this paper is to bound the degree of the 
sporadic zeros of a generic hyperplane section of the surface 
by considering the geometric implications of having sporadic 
zeros in high degree. First, we need to collect various theorems and formulae from earlier  sources.

\vspace{2mm}
1. If $S$ is a smooth surface in $\pfour$  then it satisfies 
the {\bf double point formula} ([H, pg 434]):

\begin{eqnarray}
d^2-5d-10(\pi - 1) + 2(6 \chi{\cal O}_S - K^2) & = & 0.
\end{eqnarray}

\vspace{2mm}
2.  {\bf Proposition 2. ([EP]) }

\vspace{2mm}
{\em If $S$ is a smooth surface not of general type in $\pfour$
Then for $\sigma = 5,6$ or $7$, either 
deg$\ S  \leq \frac{5(\sigma+1)(\sigma-2)}{\sigma-4}$ 
or $S$ lies on a hypersurface, $V_{\sigma}$, of degree $\sigma$.}

\vspace{2mm}
In particular, 
${\rm deg}\ S \leq 66 $  or $ S \subset V_{7}$, 
so we may assume that $S$ lies on a hypersurface of degree $7$.
Furthermore it is known ([K]),  that if $S$ lies on a hypersurface of degree
$3$ then deg$\ S \leq 8.$ 
Therefore we may assume that $S$ lies on a 
hypersurface of minimal degree $s = 4, 5, 6$ or $7$.

\vspace{2mm}
3. If $S$ is a surface not of general type and the degree of 
$S > 5$ then $K^2 \leq 9$. ([BPV])

\vspace{2mm}
4. In [BF], $\chi {\cal O}_S$ is bounded from below, using generic 
initial ideal theory, in terms of invariants arising from a generic 
hyperplane section $C$ of $S$:

\vspace{2mm}
Let ${\rm \bf C}[x_0, x_1, x_2, x_3]$ be the ring of polynomials of $\pthree$ 
under the reverse lexicographical ordering.
Let $C$ be a curve in $\pthree$, then the generic initial ideal 
of $C$,  ${\rm gin}(I_C)$, 
is generated by elements of the form $x_0^{i}x_1^{j}x_2^{k}$.

\vspace{3mm}
{\bf Definition}.
A monomial $x_0^ax_1^bx_2^c$ is a {\it sporadic zero} of $C$ 
if $x_0^ax_1^bx_2^c \notin {\rm gin}(I_C)$, but there exists 
$c' > c$ such that $x_0^ax_1^bx_2^{c'} \in {\rm gin}(I_C)$. 

\vspace{2mm}
Let $\alpha_t$ is the number of sporadic zeros in degree t and 
assume $\alpha_t = 0$ for $t > m$.

\vspace{2mm}
Let $\Gamma$ be a generic hyperplane section of $C$.
Then 
$$ {\rm gin}(I_{\Gamma}) = {\rm gin}(I_C)_{x_3 =0}^{{\rm sat}}$$
where the saturation is with respect to $x_2$. 
The generic initial ideal of $\Gamma$ is of the form 
$$ {\rm gin}(I_{\Gamma}) = (x_0^s,\  x_0^{s-1}x_1^{\lambda_{s-1}},
 \ \dots , \  x_1^{\lambda_0}), $$
where $\sum \lambda_i = d$ and 
$\lambda_{i+1}+2 \geq \lambda_i \geq \lambda_{i+1}+1$. 
The $\lm_0 > \lm_1 > \dots >  \lm_{s-1} > 0$ are called the 
{\it connected invariants}  of $\Gamma$.

\vspace{2mm}
In [BF] they show that 

\begin{eqnarray}
 \chi {\cal O}_S & \geq & 
\sum_{t=0}^{s-1} (\binom{\lm_t+t-1}{3} - \binom{t-1}{3}) 
- \sum_{t =0}^{m} \alpha_t (t-1).
\end{eqnarray} 

\vspace{3mm}
5. If $\pi$ is the genus of $C$, then 

\begin{eqnarray}
 \pi  & = & 1+ \sum_{i=0}^{s-1}(\binom{\lm_i}{2}+(i-1)\lm_i) - 
\sum_{t=0}^{m}\alpha_t .
\end{eqnarray} 

\vspace{3mm}
Combining  all these facts we obtain:

\begin{eqnarray}
 18 \geq 2K^2 & = & d^2-5d-10(\pi - 1) + 12 \chi{\cal O}_S 
  \nonumber \\
& \geq & d^2-5d-10\left(\sum_{i=0}^{s-1}(\binom{\lm_i}{2}+
(i-1)\lm_i) - \sum_{t=0}^{m}\alpha_t \right) + \nonumber \\
 & & 12\left( \sum_{t=0}^{s-1} (\binom{\lm_t+t-1}{3} - 
\binom{t-1}{3})  - \sum_{t =0}^{m} \alpha_t (t-1)\right).
\end{eqnarray} 

\vspace{2mm}
6. By the work of Gruson and Peskine ([GP]) on the numerical 
invariants of points in $\ptwo$ we have, for $d > (s-1)^2+1$

\begin{eqnarray}
1+ \sum_{i=0}^{s-1}(\binom{\lm_i}{2}+(i-1)\lm_i) & \leq & 
\frac{d^2}{2s}+(s-4)\frac{d}{2}+1 = G(d,s). 
\end{eqnarray}

7. Braun and Floystad ([BF]) show that if $s \geq 2$ and 
$d > (s-1)^2+1$

\begin{eqnarray}
\sum_{t=0}^{s-1} (\binom{\lm_t+t-1}{3} - \binom{t-1}{3}) & \geq &
s\binom{\frac{d}{s}+\frac{s-3}{2}}{3}+1-\binom{s-1}{4}.
\end{eqnarray} 

\vspace{2mm}
Thus 

\begin{eqnarray}
 18 & \geq & d^2-5d - 
10\left( \frac{d^2}{2s}+(s-4)\frac{d}{2} - 
\sum_{t=0}^{m}\alpha_t \right)  \nonumber \\
&& +12\left( s\binom{\frac{d}{s}+\frac{s-3}{2}}{3}+1-\binom{s-1}{4}
- \sum_{t =0}^{m} \alpha_t (t-1) \right)  \nonumber \\
   & = & d^2-5d-10(\frac{d^2}{2s}+(s-4)\frac{d}{2}) 
+12 s\binom{\frac{d}{s}+\frac{s-3}{2}}{3} \nonumber \\
&& + 12(1-\binom{s-1}{4}) - \sum_{t =0}^{m} \alpha_t (12t-22).  
\end{eqnarray}

\vspace{2mm}
We will use equation (7) to get an initial bound on the degree 
and then, using ${\rm Mathematica}^{\copyright}$ 
and equation (4), we will improve the bound. 

\vspace{2mm}
Thus, we need to find the smallest possible upper bound for 
$ \sum_{t =0}^{m} \alpha_t (12t-22)$ or equivalently 
$A = \sum_{t=0}^{m} \alpha_t t$. 
 
We will bound $A$ using known bounds on the number 
of sporadic zeros and by bounding the maximum degree 
of the sporadic zeros by geometric considerations. 

\vspace{2mm}
{\bf The bound on the number of sporadic zeros}

Let $\gamma = G(d,s)-\pi$.
Any bound on $\gamma$ will also bound the number of sporadic zeros 
(see  equations (3) and  (5)). 
By [EP], $\gamma \leq \frac{d(s-1)^2}{2s}$. Furthermore, 
if $S$ is a surface not of 
general type (of  degree $ > 5$), $K^2 < 6 \chi$. 
Substituting this into the double point formula $(1)$, we get  
$\pi \geq \frac{d^2-5d+10}{10}$ and thus 
\begin{eqnarray}
\sum \alpha_t & \leq & 
1 + \sum_{i=0}^{s-1}(\binom{\lm_i}{2}+(i-1)\lm_i) 
-\frac{d^2-5d+10}{10}\\
{\rm or} && \nonumber \\
\gamma  &\leq & \frac{d^2}{2s} +(s-4)\frac{d}{2} -\frac{d^2-5d}{10}.
\nonumber 
\end{eqnarray}

Taking the minimum of the bounds for $\gamma$, we get, for 
$s=4, \gamma \leq \frac{9d}{8}$, for $s=5, \gamma \leq d$,
for $s=6, \gamma \leq \frac{d(90-d)}{60}$ and for 
$s=7, \gamma \leq \frac{d(70-d)}{35}$. 

\vspace{2mm}
{\bf The bound on the degree of sporadic zeros}

For equations (7) and (4) to hold for large degree, $A$ will need to 
be  large. If there were sporadic zeros in large 
degree this would improve our chances of making $A$ large enough. 
Furthermore, every generator of ${\rm gin}(I_C)$ of the form 
$x_0^ax_1^bx_2^c$ with $c > 0$ gives us a sporadic zero in each 
degree $i$
for $a + b \leq i \leq a+ b+ c -1$. Thus we could obtain an 
upper bound on $A$ by assuming that 
there were one generator of ${\rm gin}(I_C)$ of 
the form $x_1^{\lambda_0}x_2^{z}$ where $z$ is the maximum
number of sporadic zeros. Then 
$A \leq \sum_{\lambda_0}^{\lambda_0 + z -1} t.$ But (as we saw in 
[C]) this bound is much too big.

\vspace{2mm}
Let us consider the following situation. 
Suppose $C \subset \pthree$ is a smooth curve such that 
${\rm gin}(I_C)$ has at least three generators, one of which, M, 
is of degree $r $ and all the others are 
of degree $\leq r-2$.

\vspace{2mm}
{\bf Lemma 3}

{\em C has a secant line of order r.} 

\vspace{2mm}
{\bf Proof}

Every minimal generator of ${\rm gin}(I_C)$ either arises 
from a minimal generator 
of $I_C$ or from a generator of ${\rm gin}(I_C)$  in one degree lower.  
(See [B])

Let $J$ be the ideal generated by elements of $I_C$ in degree 
$\leq r-1$. By considering the Hilbert function associated to $J$, 
we see that degree$(V(J)) = $ degree$\ (C)+1$. 
Hence, $V(J) = C \cup X$ and $X \supset L$ a line. 

Let $f$ be the generator of $I_C$ in degree r corresponding to $M$ and 
let $F = \{ f=0 \}$. By Bezout's Theorem $ F \cap L $ in $r$ points 
(up to multiplicity) and all these points must lie on $C$. Let 
$F \cap L = \sum m_ip_i $ where $p_i$ are the points of $C$. 

\vspace{2mm}
{\bf Claim}
$L$ meets $C$ at $p_i$ with multiplicity $m_i$. 

\vspace{2mm}
{\bf Proof of Claim}
$C$ is locally cut out at $p_i$ by polynomials $F_1$ and $F_2$ of 
degree $r$.  The line $L = \{l_1 = l_2 = 0\}$ meets $V(F_1)$ and 
$V(F_2)$ at $p_i$ with multiplicity $m_i$ and thus 
$$ 
{\rm length}\left(  \frac{{\cal O}_{\pthree, p_i}}{l_1, l_2, F_j} \right)  
= m_i 
$$

for $j = 1,2$. 

However 
$ \frac{{\cal O}_{\pthree, p_i}}{l_1, l_2, F_j} = 
\frac{{\cal O}_{\pone, p_i}}{F_j|_{l_1=l_2=0}} $ and 
$F_j|_{l_1 = l_2 =0} = t^{m_i}$ where $t$ is the local defining equation 
of $p_i$ in $L$. 

Hence 
$$ 
{\rm length}\left( \frac{{\cal O}_{\pthree, p_i}}{l_1, l_2, F_1, F_2} \right) = m_i 
$$

and this is the intersection multiplicity of $L$ and 
$C$ at $p_i$. \hfill{$\Box$}

\vspace{2mm}
Let us now return to the case where $S \subset \pfour$ is a smooth 
surface not of general type. 
Suppose that for generic hyperplanes $H = \{ h = 0 \}$, the 
generic hyperplane section $C_h$  of $S$ is such that 
${\rm gin}(I_{C_h})$ has at least three generators, 
one in degree $r > \frac{d}{2}$ and all others 
in degree $\leq r-2$ and hence, by Lemma 3,  $C_h$ has an 
$r$-secant line, $L_h$. 

\vspace{2mm}
{\bf Lemma 4}

{\em Generically, these secant lines are secant lines of $S$.} 

\vspace{2mm}
{\bf Proof}

Suppose for a generic $h$, $L_h \subset S$ Then for a generic 
hyperplane $H$, $S \cap H \supset L_h$. But generically 
$ S \cap H$ is a smooth irreducible curve, and hence must be 
$L_h$. But then $S = \ptwo$.    \hfill{$\Box$}

\vspace{2mm}
Thus we are in the following situation. For a generic hyperplane 
$H = \{ h =0\} $ there exists $L_h \subset H$ 
such that $L_h$ is a secant 
line of $S$ of order $r > \frac{d}{2}$.
Let $B \subset G(1,4) $ parametrize these secant lines in the 
Grassmannian of lines in $\pfour$. Let $V = \cup_{b \in B} L_b$ 
be the union of these lines in $\pfour$.

\vspace{2mm}
{\bf Proposition 5}

{\em $S$ contains a plane curve of degree $\geq r$. }

\vspace{2mm}
{\bf Proof}

As any line in $\pfour$ is contained in a 2-dimensional family 
of hyperplanes, the dimension of $B \geq 2$.

\vspace{2mm}
$V \cap S$ is at most a 2 
dimensional space, so if the dimension of $B \geq 3$, 
two lines must meet. 
Let $L_1$ and $L_2$ be  two of these intersecting lines. Let $P$ 
be the plane containing the lines. Then $S$ intersects $P$ in at 
least $2r-1 > d$ points and hence $S \cap P \supset C$ a plane 
curve.  As the secant
lines $L_i$ will meet $C$ with multiplicity at least $r$, the degree of
$C \geq r > \frac{d}{2}$. 

\vspace{2mm}
Now suppose the dimension of $B = 2$ and no two lines from $B$ 
meet. 

Let 
$$ \Phi : B \times B - - - \rightarrow {\bf {\rm P\/}^4}^* $$
send the pair $(a, b) \in B \times B$ to the hyperplane containing
$l_a$ and $l_b$. As all lines from $B$ are skew this map is well
defined away from the diagonal. 

If the dimension of the image of $\Phi$ is $4$, then there exists a generic 
hyperplane which contains two r-secant lines. However this 
contradicts the generic hyperplane section having only one r-secant 
line.

If the dimension of a  fiber is $\geq 2$, then $S$ would be 
contained in the hyperplane of the image. However $S$ is 
non-degenerate.

Thus the generic fiber is 1-dimensional and the
image of $\Phi$ is 3-dimensional. I.e. there is a 3-dimensional space of 
hyperplanes in $\pfour$ each containing a 1-dimensional family 
of skew r-secant lines of $S$. 

Let $H = \{h=0\}$ be a hyperplane in the image of $\Phi$ containing
a 1-dimensional family of skew lines. Let $S_h \subset H$ be 
the surface which is the union of these lines. 
$I_{S_h} = ( f_h, h ) $.

$$ S_h \cap S = C = \cup_{b \in B, l_b \subset H}(l_b \cap S) $$

All points of $C$ lie on an r-secant line, hence if $g \in I_C$ is a 
polynomial of degree $\leq r-1$, $g$ must vanish on $S_h$.
Therefore the only generator of $I_C$ in degree $\leq r-1$ is 
$f_h$ (and $h$). 

$S_h \subset H$ and $S \cap H = C_h$ is a hyperplane section 
of $S$ containing $C$ and hence $I_{C_h} \subset I_C$  
and so all generators of $I_{C_h}$ in degree $\leq r-1$ are
divisible by $f_h$ (mod($h$)).

Now $S$ is contained in a hypersurface  
$V_{s} \subset \pfour$ of degree $s \leq 7 (<< r-1)$, where  
$V_{s} = \{ f_{s} = 0 \}$.
Hence 
$V_{s}|_{h =0}$ contains $S_h$ for generic $H$ in the 
image of $\Phi$ and thus $V_{s}$ contains the union of all these
surfaces, which must form a three dimensional space. As 
$V_{s}$ is irreducible, any polynomial in the
ideal of $S$ of degree $\leq r-1$ must be divisible by $f_{s}$. 
Thus $I_S$ has only the one generator in degree $\leq r-1$. But
this is impossible. 

Hence there must be two secant lines meeting and as in the case of 
dim(B) = 3 or 4, we get a plane curve of degree $\geq r > \frac{d}{2}$.

\hfill{$\Box$}

\vspace{2mm}
{\bf Lemma 6}

{\em If $S$ is a smooth surface not of general type in $\pfour$ of 
degree $ d > 50$, $S$ cannot contain a plane curve of degree 
$r > \frac{d}{2}$. }

\vspace{2mm}
{\bf Proof}

Let $C \subset P$ be a plane curve  of degree 
$r > \frac{d}{2}$ contained in $S$. Let $H$ be a hyperplane containing $C$. 

Then
$S \cap H = C_h = C \cup C_{res} $. 

We have 
$$ 0 \rightarrow {\cal O}_{C \cup C_{res}} \rightarrow 
{\cal O}_{C } \oplus {\cal O}_{C_{res}} \rightarrow {\cal O}_{C \cap C_{res}} 
\rightarrow 0 $$

therefore 
$$h^1 ({\cal O}_{C_h}) \geq h^1({\cal O}_{C }) + 
h^1({\cal O}_{C _{res}})$$
and hence 
$$ g(C_h) \geq g(C) + g(C_{res}) \geq g(C). $$

$C$ is a plane curve of degree $d_C \geq \frac{d}{2}$ and so 
$$g(C) = \frac{(d_C-1)(d_C-2)}{2} - \delta 
\geq \frac{(\frac{d}{2}-1)(\frac{d}{2}-2)}{2} $$

\vspace{2mm}
On the other hand,  the Gruson-Peskine ([GP]) bound tells us that 
$$ g(C_h)  \leq \frac{d^2}{2s} + (s-4)\frac{d}{2} +1. $$

(The inequality  is true for general hyperplane sections and 
as the projective genus will stay constant it is true for all 
hyperplane sections.)

Hence
$$\frac{d^2}{2s} + (s-4)\frac{d}{2} +1 \geq 
\frac{(\frac{d}{2}-1)(\frac{d}{2}-2)}{2}. $$

This means for $s =7$, degree $\leq 42$,
for $s =6$, degree $\leq 42$ and for 
$s =5$,  degree $\leq 50.$

For $s=4 $ the inequality holds. However the Gruson-Peskine
inequality assumes there are no sporadic zeros. 
If we suppose that the number of sporadic zeros is 
$\leq \frac{3d}{4}$ then naively we have 
$A \leq \sum_{\lambda_0}^{\lambda_0 + \frac{3d}{4} -1} t.$ 
By connectedness $\lambda_0 \leq \frac{d}{4}+3$ and hence 
$A \leq \frac{5}{32}d^2 + \frac{13}{8}d-3$. Substituting back into 
equation $(7)$ we get 
$$ 0 \geq \frac{d^3}{8}-\frac{23}{8}d^2-\frac{17}{2}d+33 $$
and hence $d \leq 25$. Therefore we may assume that the 
number of sporadic zeros is $ >  \frac{3d}{4}$, then 
$ g(C_h)  < \frac{d^2}{8}  +1 -\frac{3d}{4}$ and we in fact get 
a contradiction.
\hfill{$\Box$}

\vspace{2mm}
Thus if the degree of $S > 50$, a generic hyperplane section $C$  of $S$ cannot have 
an $r$-secant line with $r > \frac{d}{2}$. In terms of the generic
initial ideal of $C$, this means that 

either 

(1) All generators of ${\rm gin}(I_C)$ are in degree $\leq \frac{d}{2}$ 

or 

(2) If there exists a generator of ${\rm gin}(I_C)$  in maximal degree 
$r > \frac{d}{2}$ then there must exist a second generator in degree $r-1$.

\vspace{2mm}
We want  to maximize $A = \sum_{t=0}^{m} \alpha_t t$ subject to 
conditions (1) and (2).
Let $z$ be the maximum number of sporadic zeros. 

$(i)$ If $\lambda_0 + z -1 \leq \frac{d}{2}$ then 
$A \leq \sum_{t= \lambda_0}^{\lambda_0 + z -1} t.$

$(ii)$ If $\lambda_0 + z -1 >  \frac{d}{2}$ 
but $\lambda_0 + \lambda_1 + z -1 \leq d $ then 
$A \leq \sum_{t= \lambda_0}^{\lfloor \frac{d}{2}\rfloor} t + 
\sum_{t= \lambda_1+1}^
{z-\lfloor\frac{d}{2}\rfloor+\lambda_0 +\lambda_1- 1} t.$

$(iii)$If $\lambda_0 + \lambda_1 + z -1 > d $ 
let $r = \lceil\frac{\lambda_0 + \lambda_1 + z}{2} \rceil$ then 
$A \leq \sum_{t= \lambda_0}^{r} t + 
\sum_{t= \lambda_1+1}^{r-1} t.$

To get a first estimate of $A$ we have, by the connectedness 
of the invariants, 
$\lambda_0 \leq \frac{d}{s} + s - 1$ and 
$\lambda_1 \leq \frac{d}{s} + s - 2$.

\vspace{2mm}
Using the bounds on $\gamma$, we get 
$$\begin{array}{ll}
{\rm for} \ s=4, & A \leq \frac{153}{256}d^2+
 \frac{45}{16}d+\frac{1}{4} \\
{\rm for} \ s=5, & A \leq \frac{9}{20}d^2+ \frac{7}{2}d+\frac{1}{4}. 
\end{array}$$

Substituting back into the original equation $(7)$ above, we get 

$$\begin{array}{lll}
{\rm for} \ s=4, & 0 \geq \frac{1}{8}d^3-\frac{523}{64}d^2-
\frac{29}{2}d- 6 & {\rm and \  hence} \ d \leq 67, \\
{\rm for} \ s=5, & 0 \geq \frac{2}{25}d^3-\frac{27}{5}d^2-32d-21 & 
{\rm and \  hence} \ d \leq 71.
\end{array}$$

\vspace{2mm}
We now need to consider equation $(4)$ which is much more 
accurate than equation $(7)$. From point 2. we  know that either
deg$\ S \leq 90$ or $S$ is contained in a hypersurface of degree $5$.
Furthermore, if $S$ is contained in a hypersurface of degree $5$ then 
$d \leq 71$ and if $S$ is contained in a hypersurface of degree $4$
then $d \leq 67$.  
Therefore we can write down all possible configurations of the connected 
invariants 
$\lambda_0 > \lambda_1 > \dots  > \lambda_{s-1}$ for high degree. 

For example if $s =5 $ and $d = 71 $ the possible invariants are
$$\begin{array}{c}
18  > 16 > 14 > 12 > 11 \\
17  > 16 > 14 > 13 > 11 \\
17 > 15 > 14 > 13 > 12. 
\end{array}$$ 

(To obtain the list we used a program of Rich Liebling.)

We then obtain an upper bound on the number of sporadic 
zeros, $z$ using 
$z \leq \frac{d^2}{8}-
\sum(\binom{\lambda_i}{2}-(i-1)\lambda_i) + \frac{9d}{8} 
$ if $s = 4$ and equation $(8)$ if $s \geq 5$.  

We again get an upper bound on $A$ using $(i)$, $(ii)$, or $(iii)$.
Substituting  everything into equation (4) and see when the 
inequality holds. 
(We checked the inequalities on ${\rm Mathematica}^{\copyright}$.)

We  get  $s \leq 7$ and 

$$\begin{array}{lll}
{\rm for} &\ s=7, & d \leq 43 \\
{\rm for} &\ s=6, & d \leq 44 \\
{\rm for} &\ s=5, & d \leq 66 \\
{\rm and \ for} & \ s = 4, & d \leq 65.
\end{array}$$

\vspace{3mm}
{\it Acknowledgements}. We would like to thank Sheldon Katz for 
inviting the first author to Oklahoma State University and Micheal 
Schneider for inviting the second author to Bayreuth University,   
enabling us to work on this project. We would also like to thank 
Gunnar Fl{\o}ystad for many helpful conversations and especially 
for pointing out the proof of the claim in Lemma 3. 

\vspace{5mm}

{\bf References} 

{\bf [B]} D. Bayer {\it The division algorithm and the Hilbert scheme}, 
Ph.D. Thesis, Harvard University (1982). 

{\bf [BF]} R. Braun, G. Fl{\o}ystad {\it A bound for the degree 
of smooth surfaces in $\pfour$ not of general type}, Compositio 
Mathematica, Vol. 93, No. 2, September(I) (1994) 211-229.

{\bf [BPV]} W. Barth, C. Peters, A. Van de Ven 
{\it Compact Complex Surfaces}, Springer-Verlag (1984). 

{\bf [C]} M. Cook {\it An improved bound for the degree of smooth 
surfaces in $\pfour$ not of general type}, to appear in Compositio
Mathematica.

{\bf [EP]} G. Ellingsrud, C. Peskine {\it Sur les surfaces lisses de 
$\pfour$}, Invent. Math. 95 (1989) 1-11.

{\bf [GP]} L. Gruson, C. Peskine {\it Genres des courbes de l'espace 
projectif}, Lecture Notes in Mathematics, Algebraic Geometry, 
Troms{\o} 1977, 687 (1977) 31-59.

{\bf [H]} R. Hartshorne {\it Algebraic Geometry}, Springer-Verlag (1977). 

{\bf [K]} L. Koelblen {\it Surfaces de $\pfour$ trac{\'e}es sur une
 hypersurfaces cubique} Journal f{\"u}r die Riene and Angewandte 
Mathematik, 433 (1992) 113-141.

\vspace{5mm}
Robert Braun \hfill\makebox[2.25in][l]{Michele Cook}

Mathematisches Institut \hfill\makebox[2.25in][l]{Department of Mathematics}

Universit{\"a}t Bayreuth\hfill\makebox[2.25in][l]{Pomona College}

D-8580 Bayreuth\hfill\makebox[2.25in][l]{610 N. College Avenue}

Germany \hfill\makebox[2.25in][l]{Claremont, CA 91711-6348}

\hfill\makebox[2.25in][l]{e-mail mcook\verb+@+pomona.edu}

\end{document}